\documentclass[12pt]{article}

	\newcommand{\mainTitle}{Some restricted sum formulas for double zeta values}

	\newcommand{\authorName}{MACHIDE, Tomoya}
	\newcommand{\organizationName}{Kinki University}
	\newcommand{\departmentName}{Interdisciplinary Graduate School of Science and Engineering}
	\newcommand{\majorName}{Research Center for Quantum Computing}
	\newcommand{\placeAddress}{3-4-1 Kowakae, Higashi-Osaka, Osaka 577-8502, Japan}
	\newcommand{\emailAddress}{E-mail: machide.t@gmail.com}
	\newcommand{\sectOne}{Introduction}
	\newcommand{\sectTwo}{Proofs}
\usepackage{geometry}               
\usepackage{ascmac}
\usepackage{amsmath}
\usepackage{amssymb}
\usepackage{amsthm}
\usepackage[hypertex, colorlinks=true, bookmarks=true]{hyperref}

	\DeclareMathOperator*{\OPlus}{\bigoplus}
	\newcommand{\nbk}[3]{#1#3#2}		
	\newcommand{\bgbk}[3]{\bigl{#1}#3\bigr{#2}}	
	\newcommand{\Bgbk}[3]{\Bigl{#1}#3\Bigr{#2}}			
	\newcommand{\bggbk}[3]{\biggl{#1}#3\biggr{#2}}			
	\newcommand{\Bggbk}[3]{\Biggl{#1}#3\Biggr{#2}}
	\newcommand{\autobk}[3]{\left#1#3\right#2}
	\newcommand{\mcbk}[4][?]{\ifx n#1\nbk{#2}{#3}{#4}\else\ifx b#1\bgbk{#2}{#3}{#4}\else\ifx B#1\Bgbk{#2}{#3}{#4}\else\ifx g#1\bggbk{#2}{#3}{#4}\else\ifx G#1\Bggbk{#2}{#3}{#4}\else\ifx a#1\autobk{#2}{#3}{#4}\else#4\fi\fi\fi\fi\fi\fi}
	\newcommand{\nsgsb}[1]{#1}		
	\newcommand{\bgsgsb}[1]{\big{#1}}	
	\newcommand{\Bgsgsb}[1]{\Big{#1}}			
	\newcommand{\bggsgsb}[1]{\bigg{#1}}			
	\newcommand{\Bggsgsb}[1]{\Bigg{#1}}
	\newcommand{\mcsgsb}[2][?]{\ifx n#1\nsgsb{#2}\else\ifx b#1\bgsgsb{#2}\else\ifx B#1\Bgsgsb{#2}\else\ifx g#1\bggsgsb{#2}\else\ifx G#1\Bggsgsb{#2}\else#2\fi\fi\fi\fi\fi}
	\newcommand{\myEqSpace}{\,}
	\newlength{\myEqSpaceLen}
	\newcommand{\bkR}[2][a]{\mcbk[#1]{(}{)}{#2}}
	\newcommand{\bkS}[2][a]{\mcbk[#1]{[}{]}{#2}}
	\newcommand{\bkB}[2][a]{\mcbk[#1]{\{}{\}}{#2}}

	\newcommand{\tmR}[4][?]{\bkR[#1]{#2} \ifx 0#4 \else \ifx 1#4 #3 \else {#3}^{#4} \fi \fi}	
	\newcommand{\tmS}[4][?]{\bkS[#1]{#2} \ifx 0#4 \else \ifx 1#4 #3 \else {#3}^{#4} \fi \fi}
	\newcommand{\tmM}[4][?]{\bkB[#1]{#2} \ifx 0#4 \else \ifx 1#4 #3 \else {#3}^{#4} \fi \fi}
		
	\newcommand{\setZ}{\mathbb{Z}} 
	\newcommand{\gpSym}[2][?]{\ifx g#1\mathfrak{S}_{#2} \else S_{#2} \fi}
	\newcommand{\gpAlt}[2][?]{\ifx g#1 \mathfrak{A}_{#2} \else A_{#2} \fi }
	\newcommand{\gpKleinF}[1][?]{\ifx g#1 \mathfrak{V} \else V \fi }

	\newcommand{\vPack}[1][10]{\vspace{-#1pt}}
	
	\newcommand{\lnA}[1][]{&  &}
	\newcommand{\lnP}[1]{\myEqSpace#1\myEqSpace}
	\newcommand{\lnAP}[2][]{& #2 &}

	\newcommand{\refH}[2]{#1\ref{#2}}
	
	\newcommand{\refEq}[1]{(\ref{#1})}

	\newcommand{\ome}{\omega}
	\newcommand{\mo}{(-1)}

	\newcommand{\opF}[3][?]{\ifx s#1#2/#3\else\ifx b#1(#2)/(#3)\else\ifx d#1\dfrac{#2}{#3}\else\frac{#2}{#3}\fi\fi\fi}

	\newcommand{\pw}[3][?]{\ifx!#3{#2}^{#3}\else#2^{#3}\fi}
	\newcommand{\id}[3][?]{#2_{#3}}
	\newcommand{\pwR}[3][a]{\ifx!#1{\bkR[#1]{#2}}^{#3}\else\bkR[#1]{#2}^{#3}\fi}
	\newcommand{\pwB}[3][a]{\ifx!#1{\bkB[#1]{#2}}^{#3}\else\bkB[#1]{#2}^{#3}\fi}
	\newcommand{\pwS}[3][a]{\ifx!#1{\bkS[#1]{#2}}^{#3}\else\bkS[#1]{#2}^{#3}\fi}
					
	\newcommand{\nFc}[3][n]{#2\bkR[#1]{#3}}
	
	\newcommand{\idFc}[4][n]{\id{#2}{#3}\bkR[#1]{#4}}

		\newcommand{\Fc}{\nFc}
	
	\newcommand{\tpT}[3][a]{ {#2}\atop \bkR[#1]{#3} }

	\newcommand{\nSm}[2][?]{\ifx l#1 \sum\limits_{#2} \else \sum_{#2} \fi}
	\newcommand{\nSmT}[3][?]{\ifx l#1 \sum\limits_{#2}^{#3} \else \sum_{#2}^{#3} \fi}	
	\newcommand{\pSm}[2][?]{\ifx t#1 \sum_{#2}^{\prime} \else \sideset{}{^\prime}\sum_{#2} \fi}
	\newcommand{\pSmT}[3][?]{\ifx t#1 \sum_{#2}^{\prime#3} \else \sideset{}{^\prime}\sum_{#2}^{#3} \fi}	
	\newcommand{\pSmN}[1][?]{\ifx t#1 \sum^{\prime} \else \sideset{}{^\prime}\sum \fi}
	\newcommand{\dSm}[2][?]{\ifx t#1 \sum_{#2}^{\dagger} \else \sideset{}{^\dagger}\sum_{#2} \fi}
	\newcommand{\dSmT}[3][?]{\ifx t#1 \sum_{#2}^{\dagger#3} \else \sideset{}{^\dagger}\sum_{#2}^{#3} \fi}	
	\newcommand{\dSmN}[1][?]{\ifx t#1 \sum^{\dagger} \else \sideset{}{^\dagger}\sum \fi}
	\newcommand{\tpTSm}[3][?]{\nSm[#1]{\tpT{#2}{#3}}}
	\newcommand{\tpTSmT}[4][?]{\nSmT[#1]{\tpT{#2}{#3} }{#4}}

		\newcommand{\Sm}{\nSm}
		\newcommand{\SmT}{\nSmT}
		\newcommand{\tpSm}{\tpTSm}
		\newcommand{\tpSmT}{\tpTSmT}
	\newcommand{\nPd}[2][?]{\ifx l#1 \prod\limits_{#2} \else \prod_{#2} \fi}
	\newcommand{\nPdT}[3][?]{\ifx l#1 \prod\limits_{#2}^{#3} \else \prod_{#2}^{#3} \fi}

	\newcommand{\nOPs}[2][?]{\ifx l#1 \OPlus\limits_{#2} \else \OPlus_{#2} \fi}
	\newcommand{\nOPsT}[3][?]{\ifx l#1 \OPlus\limits_{#2}^{#3} \else \OPlus_{#2}^{#3} \fi}	
	\newcommand{\pOPs}[2][?]{\ifx t#1 \OPlus_{#2}^{\prime} \else \sideset{}{^\prime}\OPlus_{#2} \fi}
	\newcommand{\pOPsT}[3][?]{\ifx t#1 \OPlus_{#2}^{\prime#3} \else \sideset{}{^\prime}\OPlus_{#2}^{#3} \fi}

	\newcommand{\nIs}[2][?]{\ifx l#1 \bigcap\limits_{#2} \else \bigcap_{#2} \fi}
	\newcommand{\nIsT}[3][?]{\ifx l#1 \bigcap\limits_{#2}^{#3} \else \bigcap_{#2}^{#3} \fi}	
	\newcommand{\pIs}[2][?]{\ifx t#1 \bigcap_{#2}^{\prime} \else \sideset{}{^\prime}\bigcap_{#2} \fi}
	\newcommand{\pIsT}[3][?]{\ifx t#1 \bigcap_{#2}^{\prime#3} \else \sideset{}{^\prime}\bigcap_{#2}^{#3} \fi}

	\newcommand{\nUn}[2][?]{\ifx l#1 \bigcup\limits_{#2} \else \bigcup_{#2} \fi}
	\newcommand{\nUnT}[3][?]{\ifx l#1 \bigcup\limits_{#2}^{#3} \else \bigcup_{#2}^{#3} \fi}	
	\newcommand{\pUn}[2][?]{\ifx t#1 \bigcup_{#2}^{\prime} \else \sideset{}{^\prime}\bigcup_{#2} \fi}
	\newcommand{\pUnT}[3][?]{\ifx t#1 \bigcup_{#2}^{\prime#3} \else \sideset{}{^\prime}\bigcup_{#2}^{#3} \fi}

	\newcommand{\nLm}[2][?]{\ifx l#1 \lim\limits_{#2} \else \lim_{#2} \fi}
		
	\newcommand{\glcondEnvLineHead}[1]{ \ifx*#1 \begin{eqnarray*} \else \begin{eqnarray}  \label{#1} \fi }
	\newcommand{\glcondEnvLineTail}[1]{ \ifx*#1 \end{eqnarray*} \else \end{eqnarray} \fi }
	\newcommand{\glcondDis}[1]{\ifx d#1 \displaystyle \fi}

		\newcommand{\envHLineT}[3][*]{ \glcondEnvLineHead{#1} #2&=&#3\glcondEnvLineTail{#1} }
		\newcommand{\envHLineTDef}[3][*]{ \glcondEnvLineHead{#1} #2&:=&#3\glcondEnvLineTail{#1} }
			\newcommand{\envHLine}{\envHLineT}
			\newcommand{\envHLineDef}{\envHLineTDef}
		\newcommand{\envMO}[2][*]{$\ifx d#1 \displaystyle \fi#2$}
		\newcommand{\envMT}[3][*]{$\ifx d#1 \displaystyle \fi#2=#3$}
		\newcommand{\envMTDef}[3][*]{$\ifx d#1 \displaystyle \fi#2:=#3$}
			\newcommand{\envM}{\envMT}
			
		\newcommand{\envMTh}[4][*]{$\ifx d#1 \displaystyle \fi#2=#3=#4$}
		\newcommand{\envMF}[5][*]{$\ifx d#1 \displaystyle \fi#2=#3=#4=#5$}
		
	\newcommand{\envMLineT}[3][*]{ \ifx*#1 \begin{multline*} #2\lnP{=}#3\end{multline*} \else \begin{multline} \label{#1} #2\lnP{=}#3\end{multline} \fi }
	\newcommand{\envMLineTDef}[3][*]{ \ifx*#1 \begin{multline*} #2\lnP{:=}#3\end{multline*} \else \begin{multline} \label{#1} #2\lnP{:=}#3\end{multline} \fi }

	\newcommand{\abs}[1]{\left | #1 \right |  }		
	\newcommand{\sbG}[2][n]{\bkS[#1]{#2}}
	\newcommand{\sbGf}[2][n]{\bkB[#1]{#2}}

	\newcommand{\exx}[2][n]{ \Fc[#1]{\exp}{#2} }
	\newcommand{\exN}[1]{ e^{#1}}

	\newcommand{\fcZeta}[2][n]{\Fc[#1]{\zeta}{#2}}
		\newcommand{\fcZ}{\fcZeta}		
	\newcommand{\nmB}[1]{B_{#1}}
	
	\newcommand{\envMLineTPt}[4][*]{ \ifx*#1 \begin{multline*} #3\lnP{#2}#4\end{multline*} \else \begin{multline} \label{#1} #3\lnP{#2}#4\end{multline} \fi }

	\newcommand{\cTxT}[2]{\textcolor{#1}{#2}}
	
		\newcommand{\cTx}{\cTxT}
	\newcommand{\alTx}[1]{\cTx{red}{#1}}
	\newcommand{\rmTx}[1]{\cTx{blue}{#1}}
	\newcommand{\ntTx}[1]{\cTx{Green}{#1}}
	\newcommand{\cmoTx}[1]{\cTx{Gray}{#1}}
	\newcommand{\sTx}[2][n]{ \ifx t#1{\tiny #2} \else \ifx s#1{\scriptsize #2} \else \ifx f#1{\footnotesize #2} \else \ifx S#1{\small #2} \else \ifx n#1{#2} \else \ifx l#1{\large #2} \else \ifx L#1{\Large #2} \else \ifx R#1{\LARGE #2} \else \ifx h#1{\huge #2} \else \ifx H#1{\Huge #2} \else \ifx ?#1 #2 \else #2 \fi\fi\fi\fi\fi\fi\fi\fi\fi\fi\fi }
		\newcommand{\envHLineTCl}[3][a]{ \ifx a#1\alTx{\envHLineT{#2}{#3}} \else \ifx r#1 \rmTx{\envHLineT{#2}{#3}} \else\ifx n#1 \ntTx{\envHLineT{#2}{#3}}\else\ifx c#1 \cmoTx{\envHLineT{#2}{#3}} \else \text{[argument error]} \fi\fi\fi\fi \vPack[18] }

		\newcommand{\envHLineCSClPart}[8][?]{\ifx*#1 \begin{eqnarray*} \else \begin{eqnarray}  \label{#1}  \fi \alTx{#3}&#2&\alTx{#4}\\#5&#2&#6\nonumber\\\alTx{#7}&#2&\alTx{#8}\nonumber\glcondEnvLineTail{*}}
		
			\newcommand{\HLineCTCl}[3][?]{\alTx{#2}&=&\alTx{#3}\nonumber \ifx#1p \\ \fi}
			\newcommand{\HLineCTClDef}[3][?]{\alTx{#2}&:=&\alTx{#3}\nonumber \ifx#1p \\ \fi}
			\newcommand{\HLineCFCl}[5][?]{\alTx{#2}&=&\alTx{#3}\nonumber\\#4&=&#5\nonumber \ifx#1p \\ \fi}
			\newcommand{\HLineCFClDef}[5][?]{\alTx{#2}&:=&\alTx{#3}\nonumber\\#4&:=&#5\nonumber \ifx#1p \\ \fi}
			\newcommand{\HLineCSCl}[7][?]{\alTx{#2}&=&\alTx{#3}\nonumber\\#4&=&#5\nonumber\\\alTx{#6}&=&\alTx{#7}\nonumber\ifx#1p \\ \fi}
			\newcommand{\HLineCSClDef}[7][?]{\alTx{#2}&:=&\alTx{#3}\nonumber\\#4&:=&#5\nonumber\\\alTx{#6}&:=&\alTx{#7}\nonumber\ifx#1p \\ \fi}									
			\newcommand{\HLineCECl}[9][?]{\alTx{#2}&=&\alTx{#3}\nonumber\\#4&=&#5\nonumber\\\alTx{#6}&=&\alTx{#7}\nonumber\\#8&=&#9\nonumber\ifx#1p \\ \fi}
			\newcommand{\HLineCEClDef}[9][?]{\alTx{#2}&:=&\alTx{#3}\nonumber\\#4&:=&#5\nonumber\\\alTx{#6}&:=&\alTx{#7}\nonumber\\#8&:=&#9\nonumber\ifx#1p \\ \fi}
		
	\newcommand{\envMyThm}[2][*]{\ifx*#1\begin{myThm}\else\begin{myThm}[#1]\fi #2\end{myThm}}
	\newcommand{\envMyProp}[2][*]{\ifx*#1\begin{myPorp}\else\begin{myPorp}[#1]\fi #2\end{myProp}}
	\newcommand{\envMyLem}[2][*]{\ifx*#1\begin{myLem}\else\begin{myLem}[#1]\fi #2\end{myLem}}
	\newcommand{\envMyCor}[2][*]{\ifx*#1\begin{myCor}\else\begin{myCor}[#1]\fi #2\end{myCor}}
	\newcommand{\envMyDef}[2][*]{\ifx*#1\begin{myDef}\else\begin{myDef}[#1]\fi #2\end{myDef}}

	\newcommand{\envCenter}[2][*]{\ifx*#1\begin{center}\else\begin{center}[#1]\fi #2\end{center}}
	\newcommand{\envFlushleft}[2][*]{\ifx*#1\begin{flushleft}\else\begin{flushleft}[#1]\fi #2\end{flushleft}}
	\newcommand{\envFlushright}[2][*]{\ifx*#1\begin{flushright}\else\begin{flushright}[#1]\fi #2\end{flushright}}
		
	\newcommand{\envItemIm}[2][*]{\ifx*#1\begin{itemize}\else\begin{itemize}[#1]\fi #2\end{itemize}}
	\newcommand{\envItemDp}[2][*]{\ifx*#1\begin{description}\else\begin{description}[#1]\fi #2\end{description}}
	\newcommand{\envItemEm}[2][*]{\ifx*#1\begin{enumerate}\else\begin{enumerate}[#1]\fi #2\end{enumerate}}

	\newcommand{\envMultCol}[3][*]{\begin{multicols}{#2}\ifx*#1\else\mbox{}\vspace{-#1pt}\fi#3\end{multicols}}

\theoremstyle{plain}
\newtheorem{theorem}{THEOREM}[section]
\newtheorem{proposition}[theorem]{PROPOSITION}
\newtheorem{lemma}[theorem]{LEMMA}
\newtheorem{corollary}[theorem]{COROLLARY}
\theoremstyle{definition}

\theoremstyle{remark}

\allowdisplaybreaks[4]
\numberwithin{equation}{section}

	\newcommand{\lccondBibitem}[3][]{ \if ?#2 \bibitem{#3} \else \bibitem[#2]{#3} \fi}
	\newcommand{\refPaper}[8][?]{
			\lccondBibitem{#1}{#2}%bibitem		\bibitem[#1]{#2}  		%bibitem
				#3,			%author
				\emph{#4}, 	%artName
				#5\ 			%jouName 
				{\bf #6},		%volume
				#7,			%year
				#8.			%page
		}
	
	\newcommand{\refPaperRep}[9][?]{
			\lccondBibitem{#1}{#2}%bibitem
				#3,			%author
				\emph{#4}, 	%artName
				#5\ 			%jouName 
				{\bf #6},		%volume
				#7,			%year
				#8			%page
				; reprinted in #9	%reprinted information
		}
	
	\newcommand{\refBook}[7][?]{
			\lccondBibitem{#1}{#2}%bibitem
				#3,			%author
				\emph{#4}, 	%bookName
				#5,			%publisher 
				#6,			%pubPlace
				#7.			%pubYear
		}
	\newcommand{\refPaperAlm}[5][?]{
			\lccondBibitem{#1}{#2}%bibitem
				#3,	 		%author
				\emph{#4}, 	%artName
				#5		%etc 
		}

	\newcommand{\pcstSpaceForRef}{\ }
	\newcommand{\refThm}[2][?]{\ifx?#1\refH{Theorem\pcstSpaceForRef}{#2}\else\ifx s#1\refH{Theorems\pcstSpaceForRef}{#2}\else{[argument error]}\fi\fi}
	\newcommand{\refProp}[2][?]{\ifx?#1\refH{Proposition\pcstSpaceForRef}{#2}\else\ifx s#1\refH{Propositions\pcstSpaceForRef}{#2}\else{[argument error]}\fi\fi}
	\newcommand{\refLem}[2][?]{\ifx?#1\refH{Lemma\pcstSpaceForRef}{#2}\else\ifx s#1\refH{Lemmas\pcstSpaceForRef}{#2}\else{[argument error]}\fi\fi}
	\newcommand{\refCor}[2][?]{\ifx?#1\refH{Corollary\pcstSpaceForRef}{#2}\else\ifx s#1\refH{Corollaries\pcstSpaceForRef}{#2}\else{[argument error]}\fi\fi}
	\newcommand{\refDef}[2][?]{\ifx?#1\refH{Definition\pcstSpaceForRef}{#2}\else\ifx s#1\refH{Definitions\pcstSpaceForRef}{#2}\else{[argument error]}\fi\fi}
	\newcommand{\refRem}[2][?]{\ifx?#1\refH{Remark\pcstSpaceForRef}{#2}\else\ifx s#1\refH{Remarks\pcstSpaceForRef}{#2}\else{[argument error]}\fi\fi}
	\newcommand{\refTab}[2][?]{\ifx?#1\refH{Table\pcstSpaceForRef}{#2}\else\ifx s#1\refH{Tables\pcstSpaceForRef}{#2}\else{[argument error]}\fi\fi}
	\newcommand{\rTx}[2][0.5]{ \raise#1ex\hbox{$ \displaystyle#2$} }

	\newcommand{\glcondEnvLineTailPd}[1]{.\ifx*#1 \end{eqnarray*} \else \end{eqnarray} \fi }
	\newcommand{\glcondEnvLineTailCm}[1]{,\ifx*#1 \end{eqnarray*} \else \end{eqnarray} \fi }
	\newcommand{\envProof}[2][?]{ \par\mbox{}\vspace{-5pt}\\ \ifx?#1\emph{Proof.}\else\emph{Proof of #1.}\fi \ #2 \hfill $\Box$\\ \par}
	
		\newcommand{\envLineTPd}[3][*]{ \glcondEnvLineHead{#1} & &#2\\&=&#3\nonumber \glcondEnvLineTailPd{#1} }

			\newcommand{\envLinePd}{\envLineTPd}

		\newcommand{\envLineThPd}[4][*]{ \glcondEnvLineHead{#1} & &#2\\&=&#3\nonumber \\&=&#4\nonumber \glcondEnvLineTailPd{#1} }
		\newcommand{\envLineThCm}[4][*]{ \glcondEnvLineHead{#1} & &#2\\&=&#3\nonumber \\&=&#4\nonumber \glcondEnvLineTailCm{#1} }
		\newcommand{\envHLineTPd}[3][*]{ \glcondEnvLineHead{#1} #2&=&#3\glcondEnvLineTailPd{#1} }
		
		\newcommand{\envHLineTCm}[3][*]{ \glcondEnvLineHead{#1} #2&=&#3\glcondEnvLineTailCm{#1} }
		
			\newcommand{\envHLinePd}{\envHLineTPd}
			
			\newcommand{\envHLineCm}{\envHLineTCm}
			
		\newcommand{\envHLineFPd}[5][*]{ \glcondEnvLineHead{#1} #2&=&#3\\&=&#4\nonumber \\&=&#5\nonumber \glcondEnvLineTailPd{#1} }
		\newcommand{\envHLineFCm}[5][*]{ \glcondEnvLineHead{#1} #2&=&#3\\&=&#4\nonumber \\&=&#5\nonumber \glcondEnvLineTailCm{#1} }
		\newcommand{\envHLineCFNmePd}[5][?]{\begin{eqnarray} #2&=&#3,\\#4&=&#5 \glcondEnvLineTailPd{?} }
		\newcommand{\envHLineCFDefNmePd}[5][?]{\begin{eqnarray} #2&:=&#3,\\#4&:=&#5 \glcondEnvLineTailPd{?} }
		\newcommand{\envHLineCFCmNme}[5][?]{\begin{eqnarray} #2&=&#3,\\#4&=&#5 \glcondEnvLineTailCm{?} }
		\newcommand{\envHLineCFCmDefNme}[5][?]{\begin{eqnarray} #2&:=&#3,\\#4&:=&#5 \glcondEnvLineTailCm{?} }
		\newcommand{\envHLineCSNmePd}[7][*]{\begin{eqnarray} #2&=&#3,\\#4&=&#5,\\#6&=&#7\glcondEnvLineTailPd{?}}
		\newcommand{\envHLineCSDefNmePd}[7][*]{\begin{eqnarray} #2&:=&#3,\\#4&:=&#5,\\#6&:=&#7\glcondEnvLineTailPd{?}}
		\newcommand{\envHLineCSCmNme}[7][*]{\begin{eqnarray} #2&=&#3,\\#4&=&#5,\\#6&=&#7\glcondEnvLineTailCm{?}}
		\newcommand{\envHLineCSCmDefNme}[7][*]{\begin{eqnarray} #2&:=&#3,\\#4&:=&#5,\\#6&:=&#7\glcondEnvLineTailCm{?}}
		\newcommand{\envHLineCENmePd}[9][*]{\begin{eqnarray} #2&=&#3,\\#4&=&#5,\\#6&=&#7,\\#8&=&#9\glcondEnvLineTailPd{?}}
		\newcommand{\envHLineCEDefNmePd}[9][*]{\begin{eqnarray} #2&:=&#3,\\#4&:=&#5,\\#6&:=&#7,\\#8&:=&#9\glcondEnvLineTailPd{?}}
		\newcommand{\envHLineCECmNme}[9][*]{\begin{eqnarray} #2&=&#3,\\#4&=&#5,\\#6&=&#7,\\#8&=&#9\glcondEnvLineTailCm{?}}
		\newcommand{\envHLineCECmDefNme}[9][*]{\begin{eqnarray} #2&:=&#3,\\#4&:=&#5,\\#6&:=&#7,\\#8&:=&#9\glcondEnvLineTailCm{?}}
		\newcommand{\envPLinePd}[2][*]{\glcondEnvLineHead{#1} #2\glcondEnvLineTailPd{#1}}
		\newcommand{\envPLineCm}[2][*]{\glcondEnvLineHead{#1} #2\glcondEnvLineTailCm{#1}}
		\newcommand{\envOTLinePd}[4][*]{\glcondEnvLineHead{#1} #2\lnAP{=}#3\lnP{=}#4.\glcondEnvLineTail{#1}}

		\newcommand{\envMOCm}[2][*]{$\ifx d#1 \displaystyle \fi#2$,}
		\newcommand{\envMOPd}[2][*]{$\ifx d#1 \displaystyle \fi#2$.}
		
		\newcommand{\envMTCm}[3][*]{$\ifx d#1 \displaystyle \fi#2=#3$,}
		\newcommand{\envMTPd}[3][*]{$\ifx d#1 \displaystyle \fi#2=#3$.}
		\newcommand{\envMTCmDef}[3][*]{$\ifx d#1 \displaystyle \fi#2:=#3$,}
		\newcommand{\envMTDefPd}[3][*]{$\ifx d#1 \displaystyle \fi#2:=#3$.}
			\newcommand{\envMCm}{\envMTCm}

		\newcommand{\envMThCm}[4][*]{$\ifx d#1 \displaystyle \fi#2=#3=#4$,}
		\newcommand{\envMThPd}[4][*]{$\ifx d#1 \displaystyle \fi#2=#3=#4$.}

		\newcommand{\envHLineCFCmNm}[5][*]{ \begin{equation}\begin{split} \ifx*#1 \text{[ERROR;need label name]} \else \label{#1} \fi #2&\lnP{=}#3,\\#4&\lnP{=}#5, \end{split}\end{equation} }
		\newcommand{\envHLineCFNm}[5][*]{ \begin{equation}\begin{split} \ifx*#1 \text{[ERROR;need label name]} \else \label{#1} \fi #2&\lnP{=}#3\\#4&\lnP{=}#5, \end{split}\end{equation} }
		\newcommand{\envHLineCFNmPd}[5][*]{ \begin{equation}\begin{split} \ifx*#1 \text{[ERROR;need label name]} \else \label{#1} \fi #2&\lnP{=}#3,\\#4&\lnP{=}#5. \end{split}\end{equation} }
		\newcommand{\envHLineCFCmDefNm}[5][*]{ \begin{equation}\begin{split} \ifx*#1 \text{[ERROR;need label name]} \else \label{#1} \fi #2&\lnP{:=}#3,\\#4&\lnP{:=}#5, \end{split}\end{equation} }
		\newcommand{\envHLineCFDefNm}[5][*]{ \begin{equation}\begin{split} \ifx*#1 \text{[ERROR;need label name]} \else \label{#1} \fi #2&\lnP{:=}#3\\#4&\lnP{:=}#5, \end{split}\end{equation} }
		\newcommand{\envHLineCFDefNmPd}[5][*]{ \begin{equation}\begin{split} \ifx*#1 \text{[ERROR;need label name]} \else \label{#1} \fi #2&\lnP{:=}#3,\\#4&\lnP{:=}#5. \end{split}\end{equation} }
		\newcommand{\envHLineCSCmNm}[7][*]{ \begin{equation}\begin{split} \ifx*#1 \text{[ERROR;need label name]} \else \label{#1} \fi #2&\lnP{=}#3,\\#4&\lnP{=}#5,\\#6&\lnP{=}#7 \end{split}\end{equation} }
		\newcommand{\envHLineCSNm}[7][*]{ \begin{equation}\begin{split} \ifx*#1 \text{[ERROR;need label name]} \else \label{#1} \fi #2&\lnP{=}#3\\#4&\lnP{=}#5\\#6&\lnP{=}#7 \end{split}\end{equation} }
		\newcommand{\envHLineCSNmPd}[7][*]{ \begin{equation}\begin{split} \ifx*#1 \text{[ERROR;need label name]} \else \label{#1} \fi #2&\lnP{=}#3,\\#4&\lnP{=}#5,\\#6&\lnP{=}#7. \end{split}\end{equation} }
		\newcommand{\envHLineCSCmDefNm}[7][*]{ \begin{equation}\begin{split} \ifx*#1 \text{[ERROR;need label name]} \else \label{#1} \fi #2&\lnP{:=}#3,\\#4&\lnP{:=}#5,\\#6&\lnP{:=}#7 \end{split}\end{equation} }
		\newcommand{\envHLineCSDefNm}[7][*]{ \begin{equation}\begin{split} \ifx*#1 \text{[ERROR;need label name]} \else \label{#1} \fi #2&\lnP{:=}#3\\#4&\lnP{:=}#5\\#6&\lnP{:=}#7 \end{split}\end{equation} }
		\newcommand{\envHLineCSDefNmPd}[7][*]{ \begin{equation}\begin{split} \ifx*#1 \text{[ERROR;need label name]} \else \label{#1} \fi #2&\lnP{:=}#3,\\#4&\lnP{:=}#5,\\#6&\lnP{:=}#7. \end{split}\end{equation} }
		\newcommand{\envHLineCECmNm}[9][*]{ \begin{equation}\begin{split} \ifx*#1 \text{[ERROR;need label name]} \else \label{#1} \fi #2&\lnP{=}#3,\\#4&\lnP{=}#5,\\#6&\lnP{=}#7,\\#8&\lnP{=}#9,  \end{split}\end{equation} }
		\newcommand{\envHLineCENm}[9][*]{ \begin{equation}\begin{split} \ifx*#1 \text{[ERROR;need label name]} \else \label{#1} \fi #2&\lnP{=}#3\\#4&\lnP{=}#5\\#6&\lnP{=}#7\\#8&\lnP{=}#9  \end{split}\end{equation} }
		\newcommand{\envHLineCENmPd}[9][*]{ \begin{equation}\begin{split} \ifx*#1 \text{[ERROR;need label name]} \else \label{#1} \fi #2&\lnP{=}#3,\\#4&\lnP{=}#5,\\#6&\lnP{=}#7,\\#8&\lnP{=}#9.  \end{split}\end{equation} }
		\newcommand{\envHLineCECmDefNm}[9][*]{ \begin{equation}\begin{split} \ifx*#1 \text{[ERROR;need label name]} \else \label{#1} \fi #2&\lnP{:=}#3,\\#4&\lnP{:=}#5,\\#6&\lnP{:=}#7,\\#8&\lnP{:=}#9,  \end{split}\end{equation} }
		\newcommand{\envHLineCEDefNm}[9][*]{ \begin{equation}\begin{split} \ifx*#1 \text{[ERROR;need label name]} \else \label{#1} \fi #2&\lnP{:=}#3\\#4&\lnP{:=}#5\\#6&\lnP{:=}#7\\#8&\lnP{:=}#9  \end{split}\end{equation} }
		\newcommand{\envHLineCEDefNmPd}[9][*]{ \begin{equation}\begin{split} \ifx*#1 \text{[ERROR;need label name]} \else \label{#1} \fi #2&\lnP{:=}#3,\\#4&\lnP{:=}#5,\\#6&\lnP{:=}#7,\\#8&\lnP{:=}#9.  \end{split}\end{equation} }
	\newcommand{\envMLineTPd}[3][*]{ \ifx*#1 \begin{multline*} #2\lnP{=}#3.\end{multline*} \else \begin{multline} \label{#1} #2\lnP{=}#3.\end{multline} \fi }
	\newcommand{\envMLineTCm}[3][*]{ \ifx*#1 \begin{multline*} #2\lnP{=}#3,\end{multline*} \else \begin{multline} \label{#1} #2\lnP{=}#3,\end{multline} \fi }
	\newcommand{\envMLineTDefPd}[3][*]{ \ifx*#1 \begin{multline*} #2\lnP{:=}#3.\end{multline*} \else \begin{multline} \label{#1} #2\lnP{:=}#3.\end{multline} \fi }
	\newcommand{\envMLineTCmDef}[3][*]{ \ifx*#1 \begin{multline*} #2\lnP{:=}#3,\end{multline*} \else \begin{multline} \label{#1} #2\lnP{:=}#3,\end{multline} \fi }

			\newcommand{\lccondPar}[1]{\ifx#1p \\ \fi}
			\newcommand{\OTLineCThCm}[4][?]{#2&=&#3=#4,\nonumber \ifx#1p \\ \fi}
			\newcommand{\OTLineCThPd}[4][?]{#2&=&#3=#4.\nonumber \ifx#1p \\ \fi}
			\newcommand{\OTLineCThCmDef}[4][?]{#2&:=&#3=#4,\nonumber \ifx#1p \\ \fi}
			\newcommand{\OTLineCThDefPd}[4][?]{#2&:=&#3=#4.\nonumber \ifx#1p \\ \fi}
			\newcommand{\OTLineCSCm}[7][?]{#2&=&#3=#4\nonumber#5&=&#6=#7,\nonumber \ifx#1p \\ \fi}
			\newcommand{\OTLineCSPd}[7][?]{#2&=&#3=#4\nonumber#5&=&#6=#7.\nonumber \ifx#1p \\ \fi}
	\newcommand{\envMLineTCmPt}[4][*]{ \ifx*#1 \begin{multline*} #3\lnP{#2}#4,\end{multline*} \else \begin{multline} \label{#1} #3\lnP{#2}#4,\end{multline} \fi }
	\newcommand{\envMLineTPdPt}[4][*]{ \ifx*#1 \begin{multline*} #3\lnP{#2}#4.\end{multline*} \else \begin{multline} \label{#1} #3\lnP{#2}#4.\end{multline} \fi }

	\DeclareFontEncoding{OT2}{}{}
	\DeclareFontSubstitution{OT2}{cmr}{m}{n}
	\DeclareFontFamily{OT2}{cmr}{\hyphenchar\font45}
	\DeclareFontShape{OT2}{cmr}{m}{n}{<5><6><7><8><9>gen*wncyr <10><10.95><12><14.4><17.28><20.74><24.88>wncyr10}{}
	\DeclareFontShape{OT2}{cmr}{b}{n}{<5><6><7><8><9>gen*wncyb<10><10.95><12><14.4><17.28><20.74><24.88>wncyb10}{}
	\DeclareMathAlphabet{\mathcyr}{OT2}{cmr}{m}{n}
	\DeclareMathAlphabet{\mathcyb}{OT2}{cmr}{b}{n}
	\SetMathAlphabet{\mathcyr}{bold}{OT2}{cmr}{b}{n}

	\newcommand{\gfcDZV}[3][n]{\idFc[#1]{\mathfrak{D}}{#2}{#3}}

	\newcommand{\lctO}{\ome}
	\newcommand{\loSmN}[1][?]{\sum_{\lctO}}
	\newcommand{\lEqv}[4][\,]{#2\equiv#3#1(#4)}

\allowdisplaybreaks[4]
	\title{\mainTitle}
	\author{\authorName}
	\date{}

\begin{document}
%--maketitle
\maketitle
%--abstract
\begin{abstract}
We give some restricted sum formulas for double zeta values
	whose arguments satisfy certain congruence conditions modulo $2$ or $6$,
	and also give an application to 
	identities showed by Ramanujan
	for sums of products of Bernoulli numbers with a gap of $6$.
\end{abstract}

%--body
\section{\sectOne} \label{sectOne}
The double zeta values
	are defined by  
	\envHLineDef[1_Plain_DefMZV]
	{
		\fcZ{l_1,l_2}
	}
	{
		\Sm{m_1>m_2>0} \opF{1}{ \pw{m_1}{l_1}\pw{m_2}{l_2} }
	}
	for integers $l_1\geq2, l_2\geq1$.
These values were studied in detail in \cite{GKZ06},
	and 
	interesting facts such as linear relations 
	and connections with modular forms (especially period polynomials)
	were discovered.

Historically,
	Euler \cite{Euler1775} first studied these values,
	and showed the sum formula
	\envHLinePd[1_Plain_EqSumFml]
	{
		\tpSm{l_1\geq2, l_2\geq1}{l_1+l_2=l} \fcZ{l_1,l_2}
	}
	{
		\fcZ{l}
	} 
When the weight $l=l_1+l_2$ is even,
	Gangl, Kaneko and Zagier \cite{GKZ06} 
	gave restricted analogues of the sum formula,
	more precisely, 
	proved 
	the following formulas for double zeta values with even and odd arguments.
	\envPLineCm[1_Plain_EqsRsrSumFlm]
	{
		\pSm{l_1,l_2 \equiv 0 (2)} \fcZ{l_1,l_2}
	\lnP{=}
		\opF{3}{4} \fcZ{l}
		,
		\qquad
		\pSm{ l_1,l_2 \equiv 1 (2)} \fcZ{l_1,l_2}
	\lnP{=}
		\opF{1}{4} \fcZ{l}
	}
	where $\pSm[t]{c(l_1,l_2)}$ 
	means 
	running over the integers $l_1,l_2$ satisfying 
	$l_1\geq2, l_2\geq1, l=l_1+l_2$ and the condition $c(l_1,l_2)$.
Nakamura \cite{Nakamura09} pointed out that
	the first formula of \refEq{1_Plain_EqsRsrSumFlm} yields 
	the identity showed by Euler for sums of products of Bernoulli numbers
	\envHLineCm[1_Plane_EqBerNm]
	{
		\tpSmT{j=0}{j\equiv0 (2)}{l} \binom{l}{j} \nmB{j}\nmB{l-j}
	}
	{
		- (l-1) \nmB{l}
		\qquad
		(l\geq4)
	}
	and vice versa,
	where the Bernoulli numbers $\nmB{m}$ are defined by 
	the generating function $\opF[s]{X}{(\exN{X}-1)} = \SmT{m=0}{\infty} (\opF[s]{\nmB{m}}{m!}) X^m$.

In this paper,
	we give some new restricted sum formulas for double zeta values of any weight $l$
	whose first arguments $l_1$ satisfy certain congruence conditions modulo $2$ or $6$,
	and 
	prove that an obtained restricted sum formula yields 
	identities showed by Ramanujan for sums of products of Bernoulli numbers with a gap of $6$,
	and vice versa.

The restricted sum formulas are as follows,
	which are divided into three classes according to the value of the weight modulo $3$.

%A
%:1_Thm1
\begin{theorem}\label{1_Thm1}
Let $l$ be an integer with $l\geq3$,
	and the empty sum mean $0$.
\mbox{}\\{\bf (i)} 
If $\lEqv{l}{0}{3}$,
	\envHLinePd[1_Thm1i_Eq1]
	{
		\bkR{ \pSm{l_1 \equiv 3 (6)} - \pSm{l_1 \equiv 4 (6) } - \pSm{l_1 \equiv 5 (6)}  } \fcZ{l_1,l_2} 
	}
	{
		\opF{1}{3} \, \pSm{l_1 \equiv 1 (2)} \fcZ{l_1,l_2} 
	}
{\bf (ii)} 
If $\lEqv{l}{1}{3}$,
	\envHLinePd[1_Thm1ii_Eq1]
	{
		\bkR{ \pSm{l_1 \equiv 3 (6)} + \pSm{l_1 \equiv 4 (6) } - \pSm{l_1 \equiv 5 (6)} } \fcZ{l_1,l_2}
	}
	{
		\opF{1}{3} \, \pSm{l_1 \equiv 0 (2)} \fcZ{l_1,l_2} 
	}
{\bf (iii)} 
If $\lEqv{l}{2}{3}$,
	\envHLinePd[1_Thm1iii_Eq1]
	{
		\pSm{l_1 \equiv 4 (6)} \fcZ{l_1,l_2}
	}
	{
		\opF{1}{6}\fcZ{l} - \opF{1}{3} \pSm{l_1 \equiv 1 (2)}  \fcZ{l_1,l_2} 
	}
\end{theorem}
%Z

We restate the restricted sum formulas  in the case where $l$ is even as a corollary,
	since the restated formulas include \refEq{1_Cor1iii_Eq1} which yields the identities showed by Ramanujan,
	and the other formulas seem interesting in themselves.
Restating is easily carried out by \refEq{1_Plain_EqsRsrSumFlm} and the Chinese remainder theorem.
%A
%:1_Cor1
\begin{corollary}\label{1_Cor1}
Let $l$ be an even integer with $l\geq4$.
\mbox{}\\{\bf (i)} 
If $\lEqv{l}{0}{6}$,
	\envHLinePd[1_Cor1i_Eq1]
	{
		\bkR{ \rTx{ \pSm{l_1,l_2 \equiv 3(6)} - \pSm{l_1\equiv4(6) \atop l_2\equiv2(6) } - \pSm{l_1\equiv5(6) \atop l_2\equiv1(6)} } } \fcZ{l_1,l_2} 
	}
	{
		\opF{1}{12}\fcZ{l}
	}
{\bf (ii)} 
If $\lEqv{l}{4}{6}$,
	\envHLinePd[1_Cor1ii_Eq1]
	{
		\bkR{ \rTx{ \pSm{l_1\equiv3(6) \atop l_2\equiv1(6)} + \pSm{l_1\equiv4(6) \atop l_2\equiv0(6)} - \pSm{l_1,l_2\equiv5(6)} } } \fcZ{l_1,l_2}
	}
	{
		\opF{1}{4}\fcZ{l}
	}
{\bf (iii)} 
If $\lEqv{l}{2}{6}$,
	\envHLinePd[1_Cor1iii_Eq1]
	{
		\pSm{l_1,l_2 \equiv 4 (6)} \fcZ{l_1,l_2}
	}
	{
		\opF{1}{12}\fcZ{l}
	}

\end{corollary}
%Z

Ramanujan \cite[(13)]{Ramanujan1911} (see also \cite[(13)]{Wagstaff84})
	showed the following identities for sums of products of Bernoulli numbers  with a gap of $6$,
	\envHLine[1_Pl_EqsBerFlmRama]
	{
		\tpSmT{j=0}{j \equiv m (6)}{l} \binom{l}{j} \nmB{j}\nmB{l-j}
	}
	{
		- \opF{l-1}{3} \nmB{l}
		\qquad
		(m=0,2,4)
	}
	where $\lEqv{l}{2}{6}$ and $l\geq8$.
To be precise,
	he proved only \refEq{1_Pl_EqsBerFlmRama} with $m=0$ by using identities of trigonometric functions,
	but 
	it is easily seen that 
	the three identities in \refEq{1_Pl_EqsBerFlmRama} are equivalent; 
Identities \refEq{1_Pl_EqsBerFlmRama} with $m=0$ and $2$
	are derived from the index change $j\to l-j$ each other,
	and the two identities yield \refEq{1_Pl_EqsBerFlmRama} with $m=4$ 
	and vice versa because of \refEq{1_Plane_EqBerNm}.
Note that
	Ramanujan considered Bernoulli numbers 
	to be not $\nmB{m}$ but $\abs{\nmB{m}}$ for positive even integers $m$ in \cite{Ramanujan1911},
	and that
	there is a minor misprint in \cite[(13)]{Wagstaff84},
	that is,
	the right hand side of \cite[(13)]{Wagstaff84} should be multiplied by $\nmB{6n+2}$.
Though	
	identities of Bernoulli numbers have been studied for a very long time
	and rediscovered many times,
	\refEq{1_Pl_EqsBerFlmRama} seems truly due to Ramanujan by Wagstaff's comment in  \cite[p.54]{Wagstaff84}
	(see also \cite[Chapter 5]{Berndt85} for Ramanujan's works about Bernoulli numbers).	

We have the following corollary,
	which gives a new proof of \refEq{1_Pl_EqsBerFlmRama} via double zeta values.
%A
%:1_Cor2
\begin{corollary}\label{1_Cor2}
\refEq{1_Cor1iii_Eq1} yields \refEq{1_Pl_EqsBerFlmRama} and vice versa.
\end{corollary}
%Z

In the next and final section, we prove \refThm{1_Thm1} and \refCor{1_Cor2}.

\section{\sectTwo} \label{sectTwo}
In order to prove \refThm{1_Thm1},
	we refer to
	the proof of \refEq{1_Plain_EqsRsrSumFlm} in \cite{GKZ06},
	that is,
	we will use linear combinations of special values of the polynomials which are defined by
	\envHLineDef[2_Plane_DefGfcDZV]
	{
		\gfcDZV{l}{x,y} 
	}
	{
		\pSmN x^{l_1-1}y^{l_2-1} \fcZ{l_1,l_2}
	}	
	for integers $l\geq3$.
In fact, the formulas of \refEq{1_Plain_EqsRsrSumFlm} are obtained by
	\envPLineCm
	{
		\pSm{l_1\equiv0(2)} \fcZ{l_1,l_2} 
	\lnP{=}
		\opF{\gfcDZV{l}{1,1}-\gfcDZV{l}{-1,1}}{2}
		,
		\quad
		\pSm{l_1\equiv1(2)} \fcZ{l_1,l_2} 
	\lnP{=}
		\opF{\gfcDZV{l}{1,1}+\gfcDZV{l}{-1,1}}{2}
	}
	since
	$\gfcDZV{l}{1,1}=\fcZ{l}$ due to \refEq{1_Plain_EqSumFml}
	and
	$\gfcDZV{l}{-1,1}=-\fcZ{l}/2$ if $l$ is even (see \cite[\S2]{GKZ06}).

For a real number $x$,
	let $\sbG{x}$ and $\sbGf{x}$ respectively denote the integer and fractional parts of $x$
	such that
	$x=\sbG{x}+\sbGf{x}$, $\sbG{x}\in\setZ$ and $0\leq\sbGf{x}<1$.
	
The following proposition is necessary for the proof of \refThm{1_Thm1}.
%A
%:2_Prop1
\begin{proposition}\label{2_Prop1}
For any integer $l\geq3$,
	we have
	\envHLinePd[2_Prop1_Eq1]
	{\hspace{-15pt}
		\bkR{ \rTx{ \pSm{l_1 \equiv 2l (3) \atop l_1 \equiv 1 (2)}-\pSm{l_1 \equiv 2l (3) \atop l_1 \equiv 0 (2)}-\pSm{l_1 \equiv l-1 (3)}-2\pSm{l_1 \equiv 4 (6)} } }
		\fcZ{l_1,l_2}
	}
	{
		-\sbGf[g]{\opF{l+1}{3}} \fcZ{l} + \opF{2}{3} \gfcDZV{l}{-1,1}
	}
\end{proposition}
%Z
%A
\envProof{
We see from \cite[(26)]{GKZ06} that
	\envHLinePd[2_PropP_EqAGfcDZV]
	{
		\gfcDZV{l}{x+y,y} + \gfcDZV{l}{y+x,x} 
	}
	{
		\gfcDZV{l}{x,y} + \gfcDZV{l}{y,x} + \opF{x^{l-1}-y^{l-1}}{x-y}\fcZ{l}
	}
Let $\lctO$ denote $\exx{2 \pi i /3}$.
By summing up \refEq{2_PropP_EqAGfcDZV} with $(x,y)=(1,1),(\lctO,1),(\lctO^2,1)$ 
	and by \refLem{2_Lem1} below,
	we get
	\envLinePd[2_PropP_Eq1]
	{
		\bkR{ \pSm{l_1 \equiv 1 (3)} + \pSm{l_1 \equiv 2l (3)} } \mo^{l_1-1} \fcZ{l_1,l_2}
		+
		\opF{l+1}{3}\fcZ{l}
		-
		\opF{2}{3} \gfcDZV{l}{-1,1}	
	}
	{
		\bkR{ \pSm{l_1 \equiv 1 (3)} + \pSm{l_1 \equiv l-1 (3)} }  \fcZ{l_1,l_2} 
		+ 
		\sbG[g]{\opF{l+1}{3}} \fcZ{l}
	}
A calculation shows that
	\envLineThPd[2_PropP_Eq2]
	{
		 \bkR{ \pSm{l_1 \equiv 1 (3)} + \pSm{l_1 \equiv 2l (3)} } \mo^{l_1-1} \fcZ{l_1,l_2}
		 -
		 \bkR{ \pSm{l_1 \equiv 1 (3)} + \pSm{l_1 \equiv l-1 (3)} }  \fcZ{l_1,l_2} 
	}
	{
		\bkR{\rTx{
			\pSm{l_1 \equiv 1 (3) \atop l_1 \equiv 1 (2)} - \pSm{l_1 \equiv 1 (3) \atop l_1 \equiv 0 (2)} 
			+
			\pSm{l_1 \equiv 2l (3) \atop l_1 \equiv 1 (2)}-\pSm{l_1 \equiv 2l (3) \atop l_1 \equiv 0 (2)}
			-
			\pSm{l_1 \equiv 1 (3)} - \pSm{l_1 \equiv l-1 (3)}
		}}  \fcZ{l_1,l_2}
	}
	{
		\bkR{\rTx{
			\pSm{l_1 \equiv 2l (3) \atop l_1 \equiv 1 (2)}-\pSm{l_1 \equiv 2l (3) \atop l_1 \equiv 0 (2)}-\pSm{l_1 \equiv l-1 (3)}-2 \pSm{l_1 \equiv 1 (3) \atop l_1 \equiv 0 (2)}  
		}}
		\fcZ{l_1,l_2}
	}
Since
	$\lEqv{l_1}{1}{3}$ and $\lEqv{l_1}{0}{2}$
	if and only if
	$\lEqv{l_1}{4}{6}$, 
	\refEq{2_PropP_Eq1} and \refEq{2_PropP_Eq2} prove \refEq{2_Prop1_Eq1}.
}%Z
%A
%:2_Lem1
\begin{lemma}\label{2_Lem1}
Let $\loSmN$ mean  $\Sm{x\in\bkB{1,\lctO,\lctO^2}}$.
For any integer $l\geq3$,
	we have
	\envPLinePd
	{
		\loSmN \gfcDZV{l}{x +1,1}
	\lnAP{=}
		3 \pSm{l_1 \equiv 1 (3)} \mo^{l_1-1} \fcZ{l_1,l_2}
		+
		\opF{l+1}{2} \fcZ{l}
		-
		\gfcDZV{l}{-1,1}	
	,\\
		\loSmN \gfcDZV{l}{x +1,x}
	\lnAP{=}
		3 \pSm{l_1 \equiv 2l (3)} \mo^{l_1-1} \fcZ{l_1,l_2}
		+
		\opF{l+1}{2} \fcZ{l}
		-
		\gfcDZV{l}{-1,1}	
	,\\
		\loSmN \gfcDZV{l}{x ,1}
	\lnAP{=}
		3 \pSm{l_1 \equiv 1 (3)}  \fcZ{l_1,l_2}
	,\\
		\loSmN \gfcDZV{l}{1,x}
	\lnAP{=}
		3 \pSm{l_1 \equiv l-1 (3)}  \fcZ{l_1,l_2}
	,\\
		\bkR[G]{ \loSmN \opF{ x^{l-1}-1 }{ x-1 } } \fcZ{l}
	\lnAP{=}
		3 \sbG[g]{\opF{l+1}{3}} \fcZ{l}
	}
\end{lemma}
%Z
%A
\envProof{
Let $k$ be an integer.
Because $\lctO$ is the $3$-th root of unity,
	$1 + \lctO^{k} + \lctO^{2k}$ is equal to $3$ if $\lEqv{k}{0}{3}$
	and $0$ otherwise,
	in particular,
	$1 + \lctO + \lctO^2 =0$.
By using the weighted sum formula $\pSmN[t] 2^{l_1-1} \fcZ{l_1,l_2} = \opF[s]{(l+1)\fcZ{l}}{2}$ given in \cite{OZu08},
	it follows from \refEq{2_Plane_DefGfcDZV}  that
	\envHLineFCm
	{
		\loSmN \gfcDZV{l}{x +1,1}
	}
	{
		\pSmN \bkR{ 2^{l_1-1} + (-\lctO)^{l_1-1} + (-\lctO^2)^{l_1-1} } \fcZ{l_1,l_2}
	}
	{
		\pSmN \bkR[B]{ 2^{l_1-1} - \mo^{l_1-1} + \mo^{l_1-1}  \bkR{ 1 + \lctO^{l_1-1} + \lctO^{2(l_1-1)} } } \fcZ{l_1,l_2}
	}
	{
		\opF{l+1}{2} \fcZ{l}
		-
		\gfcDZV{l}{-1,1}	
		+
		3 \pSm{l_1 \equiv 1 (3)} (-1)^{l_1-1} \fcZ{l_1,l_2}
	}
	which verifies the first equation in the lemma.
The other equations can be proved in the same way,
	and we omit the proofs.	
}%Z
We prove \refThm{1_Thm1}.
%A
\envProof[\refThm{1_Thm1}]{
We will prove only \refEq{1_Thm1i_Eq1} 
	since we can do \refEq{1_Thm1ii_Eq1} and \refEq{1_Thm1iii_Eq1} similarly.
Assume that $\lEqv{l}{0}{3}$.
Then the left hand side of \refEq{2_Prop1_Eq1} is equal to
	\envLineThCm
	{
		\bkR{ \pSm{l_1 \equiv 3 (6)}-\pSm{l_1 \equiv 0 (6) }-\pSm{l_1 \equiv 2 (6)}-\pSm{l_1 \equiv 5 (6)}-2\pSm{l_1 \equiv 4 (6) } } \fcZ{l_1,l_2}
	}
	{
		\bkR{ \pSm{l_1 \equiv 3 (6)}-\pSm{l_1 \equiv 4 (6)}-\pSm{l_1 \equiv 5 (6)}-\pSm{l_1 \equiv 0 (2) } } \fcZ{l_1,l_2}
	}
	{
		\bkR{ \pSm{l_1 \equiv 3 (6)} - \pSm{l_1 \equiv 4 (6) } - \pSm{l_1 \equiv 5 (6)} } \fcZ{l_1,l_2} 
		-
		\opF{\gfcDZV{l}{1,1}-\gfcDZV{l}{-1,1}}{2}
	}
	and the right hand side is equal to
	\envPLinePd
	{
		-\opF{1}{3}\gfcDZV{l}{1,1} + \opF{2}{3}\gfcDZV{l}{-1,1}
	}
We thus obtain
	\envHLineCm
	{
		\bkR{ \pSm{l_1 \equiv 3 (6)} - \pSm{l_1 \equiv 4 (6) } - \pSm{l_1 \equiv 5 (6)} } \fcZ{l_1,l_2} 
	}
	{
		\opF{\gfcDZV{l}{1,1}+\gfcDZV{l}{-1,1}}{6}
	}
	which proves \refEq{1_Thm1i_Eq1}.
}%Z
Finally we prove \refCor{1_Cor2}.
%A
\envProof[\refCor{1_Cor2}]{
We will derive \refEq{1_Pl_EqsBerFlmRama} from \refEq{1_Cor1iii_Eq1}.
Since the identities in \refEq{1_Pl_EqsBerFlmRama} yield each other by virtue of \refEq{1_Plane_EqBerNm},
	we may only prove \refEq{1_Pl_EqsBerFlmRama} with $m=4$.
From the harmonic relations 
	\envMCm{
		\fcZ{l_1}\fcZ{l_2}
	}{
		\fcZ{l_1,l_2} + \fcZ{l_2,l_1} + \fcZ{l}
	}
	we see that 
	\envHLineFPd
	{
		\pSm{l_1,l_2 \equiv 4 (6)} \fcZ{l_1,l_2}
	}
	{
		\opF{1}{2}\,\pSm{l_1,l_2 \equiv 4 (6)} \bkR{ \fcZ{l_1,l_2} + \fcZ{l_2,l_1} }
	}
	{
		\opF{1}{2}\,\pSm{l_1,l_2 \equiv 4 (6)} \bkR{ \fcZ{l_1}\fcZ{l_2} - \fcZ{l}}
	}
	{
		\opF{1}{2} \bkR{ \tpSmT{j=0}{j \equiv 4 (6)}{l} \fcZ{j}\fcZ{l-j} - \opF{l-2}{6} \fcZ{l} }
	}
This with \refEq{1_Cor1iii_Eq1} gives
	\envOTLinePd
	{
		\tpSmT{j=0}{j \equiv 4 (6)}{l} \fcZ{j}\fcZ{l-j} 
	}
	{
		\opF{1}{6} \fcZ{l}
		+
		\opF{l-2}{6} \fcZ{l}
	}
	{
		\opF{l-1}{6} \fcZ{l}
	}
By Euler's formula
	\envM{
		\fcZ{m}
	}{
		-\opF{(2\pi i)^m}{2m!} \nmB{m}
	}
	for any positive even integer $m$,
	we obtain \refEq{1_Pl_EqsBerFlmRama} with $m=4$.

The converse follows by the reversing the above statements.
}%Z
%--acknowledgement
\section*{Acknowledgements}
The author would like to thank Professor Masanobu Kaneko
	for introducing him the book \cite{Berndt85} and the papers \cite{Ramanujan1911,Wagstaff84}.

%--reference

%--author information
\begin{flushleft}
\mbox{}\\ \qquad \mbox{}\\ \qquad
\majorName		\mbox{}\\ \qquad
\departmentName	\mbox{}\\ \qquad
\organizationName	\mbox{}\\ \qquad
\placeAddress		\mbox{}\\ \qquad
\emailAddress
\end{flushleft}

\end{document}